\documentclass[11pt,reqno]{amsart}

\usepackage[dvipsnames]{xcolor}
\usepackage{tikz}
\usepackage{stmaryrd}
\usetikzlibrary{matrix,arrows,decorations.pathmorphing}
\usepackage{graphicx}
\usepackage{mathpazo}
\usepackage{setspace}
\usepackage{euler}
\usepackage{amssymb}
\usepackage[stable]{footmisc}
\usepackage{amsmath,mathtools}
\usepackage{fullpage}
\usepackage{caption}
\usepackage{amsthm}
\usepackage{mathrsfs}
\usepackage{thmtools}
\usepackage{blkarray}
\usepackage{multirow}
\usepackage{picinpar} 
\usepackage{tikz-cd}
\usepackage{color}
\usepackage{verbatim}
\usepackage{hyperref}
\hypersetup{colorlinks=true, linkcolor=black, citecolor = OliveGreen, urlcolor = black}
\usepackage{amssymb}
\usepackage{amsmath}
\usepackage{enumitem,colonequals}
\usepackage{color}
\usepackage{mathrsfs}
\usepackage[all]{xy}
\usepackage{comment}
\usepackage{mathpazo}
\usepackage{euler}
\usepackage[margin=1.in]{geometry}

\newcommand{\mR}{\mathbb{R}}
\newcommand{\mQ}{\mathbb{Q}}
\newcommand{\mF}{\mathbb{F}}

\newcommand{\mD}{\mathbb{D}}

\newcommand{\mP}{\mathbb{P}}

\newcommand{\mE}{\mathbb{E}}

\newcommand{\cD}{\mathcal{D}}

\newcommand{\cO}{\mathcal{O}}

\usepackage[OT2,T1]{fontenc}
\DeclareSymbolFont{cyrletters}{OT2}{wncyr}{m}{n}
\DeclareMathSymbol{\Sha}{\mathalpha}{cyrletters}{"58}

\DeclareMathSymbol{\Sha}{\mathalpha}{cyrletters}{"58}

\newcommand{\brk}[1]{ \left\lbrace #1 \right\rbrace }
\newcommand{\pwr}[1]{ \left( #1 \right) }

\newcommand{\sF}{{\mathscr F}}

\newcommand{\sL}{{\mathscr L}}
\newcommand{\sM}{{\mathscr M}}

\newcommand{\sO}{{\mathscr O}}

\RequirePackage{mathrsfs} 

\theoremstyle{theorem}
\numberwithin{equation}{subsection}
\newtheorem{thmx}{\text{Theorem}}

\newtheorem{theorem}[subsection]{Theorem}

\newtheorem{corollary}[subsection]{Corollary}
\newtheorem{conjecture}[subsection]{Conjecture}
\newtheorem{prop}[subsection]{Proposition}

\numberwithin{equation}{subsection}

\theoremstyle{definition}
\newtheorem{definition}[subsection]{\text{Definition}}
\newtheorem{example}[subsection]{Example}

\theoremstyle{remark}

\newtheorem{remark}[subsection]{Remark}

\numberwithin{equation}{subsection} \numberwithin{figure}{section}

 \DeclareMathOperator{\Spec}{Spec}

\DeclareMathOperator{\vol}{vol}

\DeclareMathOperator{\ord}{ord}

\newcommand{\PGL}{\mathrm{PGL}}

\newcommand{\cdef}[1]{\textsf{\textit{#1}}}

\newcommand\QQ{\mathbb{Q}}
\newcommand\RR{\mathbb{R}}

\renewcommand{\leq}{\leqslant}

\renewcommand{\geq}{\geqslant}

\DeclareMathOperator{\supp}{supp}
\DeclareMathOperator{\bir}{b}

\makeatletter
\@namedef{subjclassname@1991}{\emph{2020} Mathematics Subject Classification}
\makeatother

\begin{document}

\title{Connections between K-stability and Vojta's conjecture}

\author{Jackson S. Morrow}
\address{Jackson S. Morrow \\
	Department of Mathematics\\
	University of North Texas \\
	Denton, TX 76203, USA}
\email{jackson.morrow@unt.edu}

\author{Yueqiao Wu}
\address{Yueqiao Wu \\
        School of Mathematics\\
        Institute for Advanced Study\\
        Princeton, NJ 08540, USA
	}
\email{yueqiaow@ias.edu}

\begin{abstract}
In this note, we use recent advances concerning the K-stability of $\mQ$-Fano varieties to provide settings for which Vojta's conjecture holds. 
\end{abstract}


\date{\today}
\maketitle

\maketitle

\section{\bf Introduction}

In \cite{Vojta87}, Vojta described a dictionary between value distribution theory and Diophantine geometry, which culminated in several conjectures that unify many aspsects of arithmetic geometry. 
Vojta's main conjecture \cite[Conjecture 3.4.3] {Vojta87}, which we recall in \autoref{conj:Vojta}, is perhaps the deepest conjecture in arithmetic geometry. 
Loosely speaking, this conjecture is a Diophantine approximation statement, which concerns the relationship between the heights of points on varieties with respect to three divisors on a projective variety $X$ over number field $F$, namely a simple normal crossings divisor $D$, the canonical divisor $K_X$, and a big divisor $A$ on the variety. We recall the relevant definitions in the statement of Vojta's conjecture in Section \ref{sec:prelims}. 
This conjecture and its generalization \cite[Conjecture 5.2.6]{Vojta87} are known to have a range of important consequences. For example, it generalizes Faltings' theorem concerning the Mordell conjecture, and it is also known to imply the Bombieri--Lang conjecture and the abc-conjecture \cite[Section 5.5]{Vojta87}.

\subsection{Main results}
The goal of this note is to use recent breakthroughs in the subject of K-stability of $\mQ$-Fano varieties to prove new instances of Vojta's conjecture for $\mQ$-Fano varieties. 
Using the theory of local Weil functions associated to $\bir$-divisors (Section \ref{sec:birlocalWeil}), we can formulate a generalization of Vojta's conjecture (\autoref{conj:birVojta}) which incorporates $\bir$-divisors on $X$, and our first result gives examples of varieties and divisors which satisfy this conjecture. 
We refer the reader to Subsections \ref{subsec:varieties} and \ref{subsec:birdiv} for our conventions concerning varieties and $\bir$-divisors. 

\begin{thmx}\label{thmx:main1}
Let $F$ be a number field and let $X$ be a $\mQ$-Fano $F$-variety such that $X_{\overline{F}}$ has canonical singularities and infinite automorphism group. 
Then there exists a finite extension $F'/F$ and a $\bir$-divisor $\mE$ on $X_{F'}$ such that inequalities predicted by the birational Vojta's conjecture are true in the following sense. Let $\mD = \sum_{i=1}^q \mD_i$ be a $\bir$-divisor on $X_{F'}$ such that the traces of each $\mD_i$ are linearly equivalent to the trace of $\mE$ on some fixed normal proper model $Y$ of $X_{F'}$, and suppose the traces of $\mD_1,\dots,D_q$ intersect properly on $Y$. Then, the birational Vojta's conjecture is true for $(X,\mD)$ (\autoref{defn:Vojtatrue}). 
\end{thmx}

Our next results concern the classical version of Vojta's conjecture in the setting of a toric Fano variety $X$. 
For such a variety, work of Blum--Jonsson \cite{BlJ} gives a combinatorial description of the $\delta$-invariant in terms of the dot product of the barycenter of the polytope associated to $-K_X$ and the primitive generators of the one dimensional cones of $\Delta$. 
Using their work, we show that for a toric Fano variety $X$ there \textit{always} exists a torus invariant divisor $D$ for which Vojta's conjecture holds for $(X,D)$. 
More precisely, we have the following.

\begin{thmx}\label{thmx:main2}
Let $F$ be a number field and let $X$ be a toric $\mQ$-Fano projective $F$-variety. 
Then there exists a torus invariant divisor $E$ on $X$ that corresponds to a primitive generator of a one-dimensional cone of the fan associated to $X$ such that inequalities predicted Vojta's conjecture are true in the following sense.
Let $D = \sum_{i=1}^q D_i$ be a Cartier divisor on $X$ such that each $D_i$ is linearly equivalent to $E$ on $X$, and suppose that $D_1,\dots,D_q$ are in general position on $X$. Then, Vojta's conjecture is true for $(X,D)$ (\autoref{defn:Vojtatrue}). 
\end{thmx}

In \autoref{thmx:main1} and \autoref{thmx:main2}, we are only able to prove results for a single divisor (or a single $\bir$-divisor). 
The restriction to a single divisor is due to the fact that we only have an upper bound of $1$ for the $\delta$-invariant of the Fano varieties in these results. 
Our next result shows that the condition of K-semistability (\autoref{defn:Kstability}) will allow us to provide collections of torus invariant divisors on toric Fano varieties for which Vojta's conjecture is true. 

\begin{thmx}\label{thmx:main3}
Let $F$ be a number field, let $X$ be a toric $\mQ$-Fano projective $F$-variety, let $\{ D_1,\dots,D_q\}$ be a collection of torus invariant divisors on $X$ that correspond to primitive generators of one-dimensional cones of the fan associated to $X$ and lie in general position, and let $D = \sum_{i=1}^qD_i$.  If $X_{\overline{F}}$ is K-semistable, then Vojta's conjecture is true for $(X,D)$ (\autoref{defn:Vojtatrue}). 
\end{thmx}


\subsection{Ideas and relation to other works}
We now explain how results in K-stability help, and along the way we recall earlier works in a similar direction. In the Fano setting, Vojta's conjecture is more accessible due in part to recent work of Ru--Vojta \cite{RuVojta:BirationalNevanlinna}, where they showed it suffices to study the asymptotic volume constant (\autoref{defn:AVC}), or equivalently in the K-stability language, the volume of a filtration, as studied in \cite{BHJ17, BlJ, BJNA1} in much greater generality. 
More precisely, the first instance of using results of K-stability to understand the asymptotic volume constants appeared in the work of Grieve \cite{Grieve:DivisorialInstability}, where the valuative criterion for K-(un)stability (i.e., the $\delta$-invariant being less than $1$) provides a prime divisor $E$ over $X$ for which the asymptotic volume constant of $-K_X$ along $E$ is $\geq 1$. This is exactly what Ru--Vojta needed. Our main observation in this note is that the condition $\delta\leq 1$ is in fact well-understood in K-stability theory, see e.g.~\cite{BX19, LXZ} (see also \cite{HeRu:Stability} for a different upper bound, in a setting where $X$ is not necessarily Fano), and it holds for any Fano variety with infinite automorphism group.

\subsection{Organization}
In Section \ref{sec:prelims}, we establish notation and recall background on algebraic geometry and on the theory of heights. 
We state Vojta's conjecture in Section \ref{sec:Vojta}, and we also describe a slight generalization of his conjecture which incorporates local $\bir$-Weil functions. 
In Section \ref{sec:AVCArithmetic}, we review the notion of a volume of filtration on the section ring of a line bundle and relate this to the asymptotic volume constant of a big line bundle along an effective Cartier divisor.
We also describe results of Ru--Vojta and Grieve illustrating the importance of this volume in the arithmetic setting. 
Additionally, we discuss background on K-stability and recall the works that allow us to compute the asymptotic volume constant of the anti-canonical divisor of a $\mQ$-Fano variety along a $\bir$-Cartier divisor. 
Finally, in Section \ref{sec:Proofs}, we prove our results.  

\subsection{Acknowledgements}
We thank Nathan Grieve and Julie Tzu-Yueh Wang for helpful comments on a first draft.  
YW is supported by NSF Grant DMS-1926686.

\section{\bf Preliminaries}
\label{sec:prelims}

In this section, we establish conventions we use throughout the work and recall some definitions and concepts from algebraic geometry and height theory. 

\subsection{Fields}\label{subsec:fields}
We will use $F$ to denote a number field, $F'/F$ to denote a finite extension of $F$, and let $M_F$ denote the set of places of $F$. 
We also let $K$ denote an arbitrary field of characteristic zero. 
As usual, the notation $\overline{F}$ (resp.~$\overline{K}$) will refer to an algebraic closure of $F$ (resp.~$K$). \\

\subsection{Varieties and divisors}
\label{subsec:varieties}
In this work, a \cdef{$K$-variety} is a geometrically integral separated scheme of finite type over $\Spec(K)$. 
For any extension $K'/K$, we will denote the base change of a $K$-variety $X$ to $\Spec(K')$ by $X_{K'}$. 

\begin{definition}\label{defn:QFano}
    A \cdef{$\mQ$-Fano $K$-variety} $X$ is a geometrically normal, geometrically irreducible projective $K$-variety such that $X_{\overline{K}}$ has at most log-terminal singularities and such that the anti-canonical divisor $-K_X$ is an ample $\mQ$-Cartier divisor. 
\end{definition}

We will primarily consider integral divisors on normal projective $K$-varieties. We now record the definition of a collection of divisors being in general position.

\begin{definition}\label{defn:generalposition}
Let $D_1,\dots,D_q$ be effective Cartier divisors on a $K$-variety $X$ of dimension $n$. 
 \begin{enumerate}
     \item We say that $D_1,\dots,D_q$ \cdef{intersect properly} if for any subset $I \subset \brk{1,\dots,q}$ and any $x\in \bigcap_{i\in I}\supp(D_i)$, the sequence $(\phi_i)_{i\in I}$ is a regular sequence in the local ring $\sO_{X,x}$ where $\phi_i$ are the local defining equations of $D_i$ for $1\leq i \leq q$.
    \item We say that $D_1,\dots,D_q$ lie in \cdef{general position} if for any $I \subseteq \{ 1 ,\dots, q\}$, we have that 
\[
\dim\pwr{\bigcap_{i\in I}\supp(D_i)} = n - \#I
\] 
if $\# I \leq n$ and $\bigcap_{i\in I}\supp(D_i) = \emptyset$ if $\#I > n$. 
\end{enumerate}
By \cite[Remark 2.2]{RuVojta:BirationalNevanlinna}, these conditions are equivalent when $X$ is Cohen--Macaulay. 
\end{definition}

Later, we will need the notion of a local Weil function associated to $\bir$-divisor. To simplify this description, we will recall the notion of the presentation of a Cartier divisor from \cite[Section 2.2.1]{BombieriGubler:HeightsDiophantineGeometry}. 

\begin{definition}\label{defn:presentation}
    Let $X$ be a normal projective $K$-variety, let $D$ be a Cartier divisor on $X$ with associated line bundle $\sO_X(D)$ and rational section $s_D$. It is well-known that there exist base point free line bundles $\sL $ and $\sM$ on $X$ such that $\sO_X(D)\cong \sL \otimes \sM^{-1}$. For generating sections $s_0,\dots,s_n$ of $\sL$ and $t_0,\dots,t_m$ of $\sM$, we call the data
\[
\cD = (s_D; \sL, \underline{s}; \sM, \underline{t})
\]
a \cdef{presentation of $D$} where $\underline{s} = \{ s_0,\dots,s_n\}$ and $\underline{t} = \{ t_0,\dots,t_m\}$. 
\end{definition}

\subsection{Toric geometry}
\label{subsec:toric}
We will make use of the standard notation for toric varieties from \cite{Fulton}. 
In this work, a \cdef{toric variety} is a normal $K$-variety $X$ that is equipped with a faithful action of a split  algebraic torus $\mathbf{G}_{m,K}^n$ over $K$ which has a dense orbit in $X$. 
Note that by \cite[Theorem 9.2.9]{CLS} and faithfully flat descent, a toric variety is Cohen--Macaulay. 
Also, by the theory of toric varieties (see e.g., \cite{CLS, Fulton}), one can associate to a projective toric variety $X$ a complete fan $\Delta$, and we will primarily denote a toric variety by $X_{\Delta}$ where $\Delta$ is the associated fan. 




\subsection{$\bir$-divisors}\label{subsec:birdiv}
We recall the notions of models of varieties and $\bir$-Cartier divisors. 
The notion of $\bir$-divisor is originally due to Shokurov; see \cite[Definition 1.7.4 and Section 2.3]{Corti:bDiv} for details. 
The letter ``$\bir$'' refers to birational.

Let $X$ be a proper $K$-variety. 
Recall that a \cdef{model} of $X$ is a proper birational morphism $Y\to X$ defined over $K$ where $Y$ is a proper $K$-variety, and we consider the category of models of $X$. 
A \cdef{$\bir$-Cartier divisor} (resp.~\cdef{$\mQ$-$\bir$-Cartier divisor}) on $X$ is an equivalence class $\mD$ of pairs $(Y,D)$ where $Y$ is a model of $X$ and $D$ is a Cartier (resp. $\mQ$-Cartier) divisor on $Y$, and the equivalence relation is given by saying that $(Y_1,D_1)\sim (Y_2,D_2)$ if $Y_1$ dominates $Y_2$ via $\phi\colon Y_1 \to Y_2$ and $D_1 = \phi^*D_2$. We represent a $\mQ$-$\bir$-Cartier divisor $\mD$ as an equivalence class of pairs $(Y,D)$ as above. In this case, we will say that $D$ is the \cdef{trace} of $\mD$ on $Y$. We say that a $\mQ$-$\bir$-Cartier divisor $\mD$ on $X$ is \cdef{effective} if it is represented by a pair $(Y,D)$ where $D$ is an effective $\mQ$-divisor on $Y$.

\subsection{Local Weil functions and height functions}
\label{sec:localWeil}
Next, we will recall the notion of local Weil functions. 
Let $\cO_F$ denote the ring of integers of $F$, recall that $M_F$ is the set of places of $F$, and let $F_v$ denote the completion of $F$ with respect to the place $v\in M_F$.

Let $X$ be a projective $F$-variety. 
Using \cite{Vojta87}, one can associate to every Cartier divisor $D$ on $X$ and every place $v\in M_F$, a local Weil function $\lambda_{D,v}\colon X(F)\setminus \supp(D) \to \mR$ where $\supp(D)$ is the support of $D$. 
For an effective $D$, $\lambda_{D,v}$ measures the $v$-adic distance of a point to $D$. 
Using local Weil functions, we can define the height function of a divisor.

\begin{definition}
    For a divisor $D$, we define the \cdef{height function associated to $D$} to be
    \[
     h_D(x) = \sum_{v\in M_F}\lambda_{D,v}(x)
    \]
    for all $x\in X(F)\setminus \supp(D)$. 
    Up to a bounded constant function, $h_D$ is independent of the choices of local Weil functions. 
\end{definition}

\subsection{Birational local Weil functions}
\label{sec:birlocalWeil}
To conclude our preliminaries section, we describe how one can define the local Weil function of a $\bir$-Cartier divisor. 
We will follow the treatment from \cite{RuVojta:BirationalNevanlinna}. The reader may also consult \cite[Birational Weil functions]{Grieve:DivisorialInstability} for an equivalent characterization of these functions.

\begin{definition}
Let $X$ be a projective $F$-variety, and let $v\in M_F$. 
\begin{enumerate}
    \item A \cdef{local $\bir$-Weil function on $X$} (resp.~a \cdef{local $\mQ$-$\bir$-Weil function on $X$}) is an equivalence class of pairs $(U,\lambda)$ where $U$ is a non-empty Zariski open subset of $X$ and $\lambda\colon U(F_v) \to \mR$ is a function such that there exists a model $\phi\colon Y\to X$ and a Cartier divisor (resp.~$\mQ$-Cartier divisor) $D$ on $Y$ such that $\lambda\circ\phi$ extends to a Weil function for $D$ (resp.~such that $n\lambda\circ \phi$ extends to a Weil function for $nD$ for some non-zero integer $n$ for which $nD$ is a Cartier divisor). Pairs $(U,\lambda)$ and $(U',\lambda')$ are \cdef{equivalent} if $\lambda = \lambda'$ on $(U\cap U')(F_v)$. 
    \item Let $\lambda$ be a local $\bir$-Weil function on $X$, and let $\mD$ be a $\bir$-Cartier divisor on $X$. We say that $\lambda$ is a \cdef{local $\bir$-Weil function for $\mD$} if $\mD$ is represented by a pair $(Y,D)$ such that if $\phi\colon Y\to X$ is the structural morphism of $Y$, then $\lambda\circ \phi$ extends to a Weil function for $D$ on $Y$. 
\end{enumerate}

To conclude, we recall some properties of these local $\bir$-Weil functions.

\begin{prop}[\protect{\cite[Proposition 4.6]{RuVojta:BirationalNevanlinna}}]\label{prop:bWeil}
Let $X$ be a projective $F$-variety. 
\begin{enumerate}
    \item Let $\lambda$ be a local $\bir$-Weil function on $X$. Then there is a unique $\bir$-Cartier divisor $\mD$ such that $\lambda$ is a local $\bir$-Weil function for $\mD$.
    \item Let $\mD$ be a $\bir$-Cartier divisor on $X$. Then there is a local $\bir$-Weil function $\lambda$ for $\mD$. 
\end{enumerate}
\end{prop}

Using \autoref{prop:bWeil}, we say that $\lambda_{\mD,v}(\cdot)$ is the local $\bir$-Weil function for a $\bir$-Cartier divisor $\mD$. 
In different settings, we will consider $\lambda_{\mD,v}$ as a function defined on a dense Zariski open subset of $X$ or on a dense Zariski open subset of $Y$ where $(Y,D)$ is a presentation of $\mD$. It will be clear from context what we are considering as the domain of $\lambda_{\mD,v}$. 
\end{definition}

To conclude this section, we describe an important example of a local $\bir$-Weil function associated a $\bir$-divisor $\mD$ on $X$ defined over a finite extension extension $F'/F$. 

\begin{example}\label{example:localbir}
    Let $X$ be a projective $F$-variety, and let $\mD$ be a $\bir$-divisor on $X$ defined over a finite extension $F'/F$. 
    Suppose that $\mD$ is represented by the pair $(Y,D)$ where $Y \to X_{F'}$ is a model of $X_{F'}$ and $D$ is a Cartier divisor on $Y$, so $D$ is the trace of $\mD$ on $Y$. 
    Let $\cD = (s_D; \sL, \underline{s}; \sM, \underline{t})$ be a presentation for the divisor $\sO_Y(D)$ for $Y$ a model of $X$ defined over $F'$. 
    For any place $w\in M_{F'}$ such that $w\mid v$, the local $\bir$-Weil function associated to $\mD$ is defined for $x \in Y\setminus \supp(D)$ as
    \[
    \lambda_{\mD,v}(x) = \max_{k}\min_{\ell} \log\left| \frac{s_k}{t_{\ell}s_D}(x)\right|_w = \max_{k}\min_{\ell} \log\left| N_{F'_w/F_v}\left(\frac{s_k}{t_{\ell}s_D}(x)\right)\right|_v^{1/[F'_w:F_v]}
    \]
    where $N_{F'_w/F_v}$ is the norm map from $F'_w\to F_v$ and $|\cdot|_v$ is the standard absolute value associated to $v$. 
    Note that when $F' = F$ and $\mD$ is a Cartier divisor on $X$, we recover the construction of a local Weil function associated to $D$. 
\end{example}

\section{\bf Vojta's conjecture}
\label{sec:Vojta}
In this section, we recall Vojta's conjecture \cite{Vojta87} and state a generalization of this conjecture which incorporates $\bir$-divisors. First, we start with Vojta's conjecture, which is one of the deepest conjectures in Diophantine geometry. 

\begin{conjecture}[Vojta's Conjecture {\cite[Conjecture 3.4.3]{Vojta87}}]\label{conj:Vojta}
Let $F$ be a number field, 
let $X$ be a smooth proper $F$-variety with canonical divisor $K_X$, 
let $D$ be a normal crossings divisor on $X$ defined over $F$, 
let $S$ be a finite set of places of $k$, 
and let $A$ be a big divisor class on $X$. 
For any $\varepsilon>0$, there exists a Zariski-closed subset $Z = Z(F,X,D,\varepsilon,S,A)$ of $X$ such that 
\begin{equation}\label{eqn:Vojta}
\sum_{v\in S}\lambda_{D,v}(x) + h_{K_X}(x) \leq \varepsilon h_A(x) + O(1)    
\end{equation}
for all points $x \in X(F)\setminus Z(F)$. 
\end{conjecture}

To conclude this section, we present a ``birational'' version of Vojta's conjecture. By ``birational'', we mean that the proximity function with respect to $D$ and $S$ from \autoref{conj:Vojta} is replaced by the the proximity function with respect to $\mD$ and $S$ where $\mD$ is some  $\bir$-Cartier divisor over $X$. 
Since local Weil functions are local $\bir$-Weil functions, this ''birational'' version of Vojta's conjecture is a slight generalization of \autoref{conj:Vojta}. 

\begin{conjecture}[Birational Vojta's Conjecture ]\label{conj:birVojta}
Let $F$ be a number field, 
let $X$ be a smooth proper $F$-variety with canonical divisor $K_X$, 
let $\mD$ be an effective Cartier $\bir$-divisor on $X$ defined over a finite extension $F'/F$ such that if $(Y,D)$ represents $\mD$, then $D$ is a simple normal crossings divisor on $Y$, 
let $S$ be a finite set of places of $F$, 
and let $A$ be a big divisor class on $X$. 
For any $\varepsilon>0$, there exists a Zariski-closed subset $Z = Z(F,X,\mD,\varepsilon,S,A)$ of $X$ such that 
\begin{equation}\label{eqn:birVojta}
\sum_{v\in S}\lambda_{\mD,v}(x) + h_{K_X}(x) \leq \varepsilon h_A(x) + O(1)   
\end{equation}
for all points $x \in X(F)\setminus Z(F)$ where $\lambda_{\mD,v}$ is the local $\bir$-Weil function defined in \autoref{example:localbir}.
\end{conjecture}

\begin{remark}
We note that our ``birational'' Vojta's conjecture does not concern the behavior of \autoref{conj:Vojta} under birational modificiations of $X$ as the height functions $h_{K_X}$ and $h_A$ are defined for divisors $K_X$ and $A$ on $X$ and not on some model of $X$. 
For a discussion on this topic, we refer the reader to \cite[Example 3.5.4]{Vojta87}. 
\end{remark}

 In order to simplify certain statements, we will introduce a definition which capture when the inequality in \autoref{conj:Vojta} and \autoref{conj:birVojta} hold.

 \begin{definition}\label{defn:Vojtatrue}
     Let $F$ be a number field, let $X$ be a normal projective $F$-variety, and let $D$ (resp.~$\mD$) be an effective Cartier divisor (resp.~an effective Cartier $\bir$-divisor) on $X$. 
     If the inequality \eqref{eqn:Vojta} (resp.~\eqref{eqn:birVojta}) holds for $F,X,D$ (resp.~$F,X,\mD$) and any choices of $\varepsilon, S,A$ outside of a Zariski-closed subset $Z(F,X,D,\varepsilon,S,A)$ (resp.~$Z(F,X,\mD,\varepsilon,S,A)$), then we say that \cdef{Vojta's conjecture} (resp.~\cdef{the birational Vojta's conjecture}) \cdef{holds for $(X,D)$} (resp.~\cdef{$(X,\mD)$}). 
 \end{definition}

\section{\bf Volume of a filtration}
\label{sec:AVCArithmetic}
In this section, we begin by recalling the notion of the volume of a filtration, see for instance \cite{BHJ17, BlJ, BJNA1} for more details. Special attention will be paid to filtrations arising from an effective Cartier divisor or a prime divisor over $X$. Then we state what is needed and known about the volume of a filtration in these two special cases.
Unless otherwise stated, $X$ is a normal projective $\overline{K}$-variety, and $\sL$ is a big line bundle on $X$.
Also, we will freely use results concerning the volume of a line bundle from \cite[Section 2.2.C]{Lazarsfeld:Positivity1}, and when dealing with volumes of divisors, we will need to work with $\mR$-divisors as discussed in \cite[Section 1.3]{Lazarsfeld:Positivity1}.

\subsection{General definitions and arithmetic importance}
Denote by $R=\bigoplus_m R_m=\bigoplus_mH^0(X, m\sL)$ the section ring of $\sL$. By a \cdef{filtration} $\sF$ on $R$ we mean a decreasing, multiplicative filtration on $R$ satisfyting $\sF^0R=R$. To such a filtration and $t\in \RR_{\geq 0}$, set
\[
    \vol(\sF^{(t)}R)\coloneqq \lim_{m\to \infty} \frac{\dim \sF^{mt}R_m}{m^n/n!}.
 \]
The filtration $\sF$ is \cdef{linearly bounded} if $\sup\{t: \vol(\sF^{(t)}R)>0\}<\infty$. 

\begin{definition}\label{defn:volumenfiltration}
To each linearly bounded filtration, one can define the \cdef{volume} of $\sF$ to be 
\[
\beta(\sL, \sF)\coloneqq\frac{1}{\vol(\sL)} \int_0^\infty \vol(\sF^{(t)}R) dt.
\]    
\end{definition}

The following two examples show the two kinds of filtration that we will consider.

\begin{example}[Asymptotic volume constants]\label{defn:AVC}
    Let $D$ be an an effective Cartier divisor on $X$. As first appeared in the works \cite{Autissier:GeometryPointsEntire} and \cite{McKinnonRoth:SeshadriConstantsDiophantineApprox}, the \cdef{asymptotic volume constant of $\sL$ along $D$}, which is denoted by $\beta(\sL, D)$ in \cite{RuVojta:BirationalNevanlinna}, is the volume of the filtration 
    \[
    \sF^mH^0(X, \sL^N)\coloneqq H^0(X, \sL^N-mD).
    \]
    It is not hard to see that $\beta(\sL, D)=\liminf_{N\to\infty} \frac{\sum_{m\geq 1} h^0(X,\sL^N(-mD))}{Nh^0(X,\sL^N)}$, which is the definition in \cite{RuVojta:BirationalNevanlinna}. See \cite[Corollary 2.12]{BlJ} for a proof. 
\end{example}

\begin{example}[$S$-invariant in K-stability]\label{defn:betaKStab}
    Let $E\subset Y\to X$ be a prime divisor over $X$ where $Y$ is a normal model of $X$. Then the associated divisorial valuation $\ord_E$ defines a filtration via
    \[\sF_E^\lambda R_m\coloneqq \{s\in R_m: \ord_E(s)\geq \lambda\}.\]
    We write $\beta(\sL, E)\coloneqq \beta(\sL, \sF_E)$. 
    In the case when $X$ is $\QQ$-Fano, we put $\beta(-K_X, E)\coloneqq \frac 1r \beta(\sL, E)$, where $r$ is such that $\sL=-rK_X$ is Cartier. It is also worth noting that $\beta$ is pullback invariant under birational morphisms.
\end{example}

For our purposes, it is important to know when $\beta\geq 1$, as suggested by the following result.

\begin{theorem}[\protect{\cite[Corollary 1.13]{RuVojta:BirationalNevanlinna} \& \cite[Theorem 4.1]{Grieve:DivisorialInstability}}]\label{thm:GrievebirVojta}
    Let $F$ be a number field, and let $S$ be a finite set of places of $F$ containing all the Archimedean places. 
    Let $X$ be a $\mQ$-Fano $F$-variety, and let $D_1,\dots,D_q$ be reduced irreducible Cartier divisors over $X$ which are defined over a finite extension $F'/F$. 
    Let $\mD_1,\dots,\mD_q$ be the $\bir$-Cartier divisor determined by $D_1,\dots,D_q$, respectively, let $\phi\colon Y \to X_{F'}$ be a normal proper model such that $D_1,\dots,D_q$ is the trace of $\mD_1,\dots,\mD_q$, and let $\mD = \mD_1 + \cdots + \mD_q$ and $D = D_1 + \cdots + D_q$. 
    If $D_1,\dots,D_q$ intersect properly (\autoref{defn:generalposition}) and $\beta(-K_{X_{F'}},D) \geq 1$, then the birational version of Vojta's conjecture (\autoref{conj:birVojta}) is true for $(X,\mD)$ i.e., for any $\varepsilon>0$, there exists a Zariski-closed subset $Z = Z(F,X,\mD,\varepsilon,S,A)$ of $X$ such that 
\[
\sum_{v\in S}\lambda_{\mD,v}(x) + h_{K_X}(x) \leq \varepsilon h_A(x) + O(1)
\]
for all points $x \in X(F)\setminus Z(F)$ where $\lambda_{\mD,v}$ is the local $\bir$-Weil function defined in \autoref{example:localbir}. 
\end{theorem}

\begin{remark}
Note that when $D_1,\dots,D_q$ are reduced irreducible Cartier divisor on $X$ and $F = F'$, then \autoref{thm:GrievebirVojta} is equivalent to \cite[Corollary 1.13]{RuVojta:BirationalNevanlinna}. 
\end{remark}

\subsection{K-stability input}
We now summarize recent results in K-stability which give $\beta\geq 1$. From now on, we assume that $X$ is $\QQ$-Fano variety (\autoref{defn:QFano}).

The volume of a prime divisor over $X$ in \autoref{defn:betaKStab} is related to K-stability of Fano varieties via the following valuative criterion. 
\begin{definition}\label{defn:delta}
    Let $X$ be a $\QQ$-Fano $\overline{K}$-variety.
    The \cdef{$\delta$-invariant of $X$} is defined to be 
    \[\delta(X)=\delta(X, -K_X)\coloneqq \inf_E\frac{A_X(E)}{\beta(-K_X, E)},\]
    where the infimum runs through all prime divisors over $X$, and $A_X(E)=1+\ord_E(K_{Y/X}).$
\end{definition}

\begin{definition}[\cite{Fujita:ValuativeCriteron, Li:Ksemistability, BlJ, LXZ}]\label{defn:Kstability}
    A $\QQ$-Fano $\overline{K}$-variety $X$ is K-\cdef{semistable} (resp. K-\cdef{stable}) if $\delta(X)\geq1$ (resp. $\delta(X)>1$), and K-\cdef{unstable} otherwise.
\end{definition}

\begin{theorem}\label{thm:QFanoinfinitedelta}
    Let $X$ be a $\QQ$-Fano $\overline{K}$-variety with infinite automorphism group. Then $\delta(X)\leq 1$, and there is a prime divisor $E$ over $X$ computing $\delta(X)$ i.e., 
    \[
    \delta(X) = \frac{A_X(E)}{\beta(-K_X,E)}.
    \]
\end{theorem}
\begin{proof}
    The result should be well-known to experts in K-stability. 
    We will recall the proof for the readers who are not familiar with the K-stability literature. 
    If $\delta(X)>1$, then by \cite[Corollary 1.3]{BX19} the automorphism group of $X$ has to be finite, contradicting our assumption. The equality in the latter claim is a consequence of \cite[Theorem 1.2]{LXZ} (see also \cite{XZ22} for a different proof).
\end{proof}

\begin{corollary}\label{coro:divisorialbetage1}
    Let $X$ be a $\QQ$-Fano $\overline{K}$-variety with infinite automorphism group and canonical singularities. 
    Then there is some divisorial valuation $v = \ord_E$ with $\beta(-K_X,E)\geq 1.$
\end{corollary}
\begin{proof}
    By \autoref{thm:QFanoinfinitedelta}, we may choose a prime divisor $E$ over $X$ with $A_X(E)\leq \beta(-K_X, E)$.
    Since $X$ has canonical singularities, we have that $A_X(E)\geq 1$, and hence $\beta(-K_X,E)\geq 1$.
\end{proof}

\subsection{The toric case}
To conclude, we recall work of Blum--Jonsson \cite{BlJ} which provides a combinatorial description of the $\delta$-invariant in the setting where $X$ is a toric $\mQ$-Fano $F$-variety. 
Recall our conventions for toric varieties (cf.~Subsection \ref{subsec:toric}). 
Fix a projective toric $F$-variety $X = X_\Delta$ associated to the fan $\Delta$ and suppose that $K_X$ is $\mQ$-Cartier. 
Let $v_1, \dots, v_d$ be primitive generators of one-dimensional cones of $\Delta$, and let $D_1, \dots, D_d$ be the corresponding toric invariant divisors. Then $-K_X=D_1+\dots + D_d$. 
The polytope associated to $-K_X$ is given by 
\[P_{-K_X} = \left\{ u\in M_\RR: \langle u, v_i\rangle \geq -1, i=1, \dots, d \right\}.\]

Under this setup, we recall the following results from \cite{BlJ}.

\begin{prop}[{\cite[Corollary 7.7 \& 7.16]{BlJ}}]\label{prop:delta}
    With notation as above,
    \[\delta(X) = \min_{i = 1,\dots,d}\frac{1}{\beta(-K_X, D_i)} = \min_{i = 1,\dots,d} \frac{1}{\langle \bar u, v_i\rangle +1},\]
    where $\overline{u} $ is the barycenter of $P_{-K_X}$.
\end{prop}

\begin{corollary}[{\cite[Corollary 7.17, 7.19]{BlJ}}]\label{coro:BlumJonssonbeta}
    If $X$ is a toric $\mQ$-Fano variety, then $\delta(-K_X)\leq 1$. The equality is achieved of and only if $X_{\overline{F}}$ is K-semistable or equivalently if the barycenter of the polytope associated to $-K_X$ is equal to the origin, in which case $\beta(-K_X, D_i)=1$ for all $i$.
\end{corollary}
\begin{proof}
    We apply \autoref{prop:delta} with $b_i=1$ for all $i$. Then there is some $i$ such that $\langle \overline{u}, v_i\rangle\geq 0$. Indeed, if $\overline{u}\neq 0$, then there must be some $i$ for which $\langle \overline{u}, v_i\rangle >0$, since otherwise all $v_i$ would lie in some half-space, contradicting the fan being complete; if $\overline{u}=0$, then $\beta(-K_X, D_i)=1$ for all $i$.
\end{proof}

\begin{remark}
    Note that $\Delta$ being complete implies that we cannot have $\beta(-K_X, D_i)>1$ for all $i$.
\end{remark}

\section{\bf Proof of Main Theorems and Examples}
\label{sec:Proofs}
In this section, we prove our main theorems and provide examples of varieties satisfying the conditions of our main theorems. 

\begin{proof}[Proof of \autoref{thmx:main1}]
    This follows from \autoref{coro:divisorialbetage1} and \autoref{thm:GrievebirVojta}. 
\end{proof}

\begin{proof}[Proof of \autoref{thmx:main2}]
    This follows from \autoref{coro:BlumJonssonbeta} and \autoref{thm:GrievebirVojta}. 
\end{proof}

\begin{proof}[Proof of \autoref{thmx:main3}]
    This follows from \autoref{coro:BlumJonssonbeta} and \autoref{thm:GrievebirVojta}.  
\end{proof}

 To conclude, we give examples of or provide references to literature describing varieties satisfying the assumptions of these results.

 \subsection{$\mQ$-Fano varieties with canonical singularities and infinite automorhpism group}
 These varieties satisfy the condition of \autoref{thmx:main1}.
 In dimension 1, the situation is not very interesting as the only $\mQ$-Fano curve is $\mP^1$ which has automorhpism group $\PGL_2$, which is infinite. 
 In dimension 2, the work \cite{CheltsovProkhorov:DelPezzoInfinite} gives a full classification of del Pezzo surfaces with Du Val singularities (or equivalently canonical singularities) that have an infinite automorphism group. 
 In dimension 3, the works \cite{Kuznetsovetal:Fano3FoldInfiniteAut, Fano3Infinite} provide a classification of smooth $\mQ$-Fano threefolds which have infinite automorhpism group. 

 \subsection{Toric $\mQ$-Fano varieties}
 These varieties satisfy the condition of \autoref{thmx:main2}.
 The dimension 1 situation is the same as above. 
 In dimension 2, the smooth toric $\mQ$-Fano varieties are the well-known smooth toric del Pezzo surfaces:~$\mP^2$, $\mP^1 \times \mP^1$, the rational normal scroll $\mF_1$, and $\mP^2$ blown-up at two and three points. 
 In dimension 3, Batryev \cite{Batyrev:ToricFano3} and Wantanabe--Wantanabe \cite{Wantanabe:ToricFano3} have completely classified the smooth toric $\mQ$-Fano 3-folds, and also in dimension 4, Batryev \cite{Batyrev:ToricFano4Fold} classified smooth toric $\mQ$-Fano 4-folds. 
 A classification in higher dimensions is unknown, but we refer the reader to \cite{Sato:FanoToricHigher} for a discussion on this problem as well as a method for characterizing toric $\mQ$-Fano varieties. 

 \subsection{K-semistable toric $\mQ$-Fano varieties}
 These varieties satisfy the conditions of \autoref{thmx:main3}. 
 Again, the dimension 1 setting is the same as above. 
 In \cite{Araujoetal:Calabi}, the authors studied the Calabi problem, which aims to classify Fano 3-folds which are K-polystable. 
 In doing so, the authors gives a complete classification of smooth toric $\mQ$-Fano surfaces and 3-folds that are K-semistable. 

 \begin{prop}[\protect{\cite[Section 2]{Araujoetal:Calabi}}]
     Let $X$ be a smooth toric $\mQ$-Fano projective $F$-surface. 
     $X_{\overline{F}}$ is K-semistable if and only if $X_{\overline{F}}$ is isomorphic to $\mP^2$, $\mP^1 \times \mP^1$, or the blow-up of $\mP^2$ at three points. 
 \end{prop}

 \begin{theorem}[\protect{\cite[Table 3.1]{Araujoetal:Calabi}}]
     Let $X$ be a smooth toric $\mQ$-Fano projective $F$-threefold. 
     $X_{\overline{F}}$ is K-semistable if and only if $X_{\overline{F}}$ is isomorphic to $\mP^3$, $\mP^1 \times \mP^2$, the blow-up of $\mP^3$ are two disjoint lines, $\mP^1 \times \mP^1 \times \mP^1$, or $\mP^1 \times S$ where $S$ is a smooth toric del Pezzo surface with $K_S^2 = 6$. 
 \end{theorem}

  \bibliography{refs}{}
\bibliographystyle{amsalpha}
\end{document}